\documentstyle[12pt,twoside]{article}
\newcommand{\rsp}{{\bf R}}
\newcommand{\csp}{{\bf C}}
\newcommand{\nsp}{{\bf N}}

\newcommand{\no}[1]{(#1)}

\newcommand{\theorem}[1]{\bigskip\par{\bf Theorem #1.}\it}
\newcommand{\lemma}[1]{\bigskip\par{\bf Lemma #1.}\it}
\newcommand{\remark}[1]{\bigskip\par{\bf Remark #1.}}

\def\endproc{\medskip\rm}
\newcommand{\proof}{\bigskip\par{\bf Proof. }}
\def\endproof{\fbox{\rule{0mm}{2mm}\hspace{3mm}} \bigskip}

\newcommand{\refs}{\bigskip\par\noindent\footnotesize\bf
REFERENCES\rm\par
\begin{enumerate}}
\def\endrefs{\end{enumerate}\normalsize}

\renewcommand{\no}{{\nonumber}}
\newcommand{\ov}{\overline}
\newcommand{\bear}{\begin{array}}
\newcommand{\enar}{\end{array}}
\newcommand{\beq}{\begin{equation}}
\newcommand{\eeq}{\end{equation}}
\newcommand{\beqn}{\begin{eqnarray}}
\newcommand{\eeqn}{\end{eqnarray}}
\newcommand{\beit}{\begin{itemize}}
\newcommand{\eeit}{\end{itemize}}

\renewcommand{\l}{\lambda}
\newcommand{\dl}{\delta}

\newcommand{\s}{\sigma}

\newcommand{\Om}{\Omega}
\renewcommand{\phi}{\varphi}

\renewcommand{\L}{\Lambda}

\newcommand{\what}{\widehat}

\newcommand{\wtil}[1]{\widetilde{#1}}
\newcommand{\bs}{\backslash}
\newcommand{\pn}{\par \noindent}
\newcommand{\med}{\medskip}
\newcommand{\qq}{\qquad}
\newcommand{\q}{\quad}

\newcommand{\media}[1]{\kern 0.4ex {-} \kern -2.0 ex {\int}_{\kern -0.9 ex {#1}}\;}

\title{Some identification problems\\ for integro-differential operator equations
\footnote{Work partially supported by the Italian Ministero dell'Universit\`a e
della Ricerca Scientifica e Tecnologica.} }

\author{ Alfredo Lorenzi (Milan)\footnote{The first author is a member of GNAFA
of the Italian Istituto Nazionale Di Alta Matematica (INDAM).}, Alexander Ramm
(Manhattan)}
\date{}

\begin{document}
\maketitle
\pn
{\bf Abstract.} We consider, in a Hilbert space $H$, the convolution
integro-differential equation $u''(t)-h*Au(t)=f(t)$, $0\le t\le T$,
$h*v(t)=\int_0^t h(t-s)v(s)\,ds$, where $A$ is a linear closed densely defined
(possibly selfadjoint and/or positive definite) operator in $H$. Under suitable
assumptions on the data we solve the inverse problem consisting of finding the
kernel $h$ from the extra data (measured data) of the type $g(t):=(u(t),\varphi)$,
where $\varphi$ is some eigenvector of $A^*$. An inverse problem for the
first-order equation $u'(t)-l*Au(t)=f(t)$, $0\le t\le T$, is also studied when
$A$ enjoys the same properties as in the previous case.
\med
\pn
\med
\pn

\section{Introduction}
\setcounter{equation}{0}
Let $A:{\cal D}(A)\to H$ be a closed linear operator densely defined in a
Hilbert space $H$ with scalar product $(\cdot,\cdot)$, norm $\|\cdot\|$, and
let ${\cal D}(A)$ be the  domain of definition of $A$. Let us assume that
there exists $\phi \in {\cal D}(A^*)$ such that
\beqn
%(1.1)
A^*\phi = \l_0\phi,
\eeqn
\beqn
%(1.2)
\l_0\in \rsp \backslash \{0\},
\eeqn
\beqn
%(1.3)
\|\phi\|=1.
\eeqn
\remark{1.1} One may assume that $\l_0\in {\mathbf C }\backslash \{0\}$:
our method remains valid in this case as well. If $\l_0$
is complex, then one has to replace it by $\overline {\l_0}$
in several formulas, such as (2.1), (2.2), (2.4), (2.5).

Assumptions (1.1)--(1.3) are satisfied if $A$ is selfadjoint, $A=A^*$,
and $A$ has an eigenvalue $\l_0\neq 0$. For example, this property holds if
$A$ is an elliptic semibounded from below selfadjoint operator in a bounded
(smooth) domain. In this case the spectrum of $A$ is discrete and consists
of a sequence of real eigenvalues $\{ \l_n(A)\}_{n=1}^{+\infty}$ going to
$+\infty$.
\med
\pn
Consider the direct problem
\beqn
%(1.4)
&& u''(t) = \int_0^t h(t-s)Au(s)\,ds + f(t),\qq 0\le t\le T,\\[2mm]
%(1.5)
&& u(0)=u_0,\q u'(0)=u_1
\eeqn
where $T>0$, and
\beqn
%(1.6)
&& h\in C([0,T])=C([0,T];\rsp),\q  f\in C^1([0,T];H),\\[2mm]
%(1.7)
&& u_0\in {\cal D}(A),\q u_1\in H,\q (u_0,\phi)\neq 0.
\eeqn
We emphasize that in equation (1.6) operator $A$ appears only
under the integral
sign. In other words, we are concerned with the case when the differential
operator outside the integral does {\it not} dominate the operator inside
the integral.
\pn
Assume now that the kernel $h$ is a real-valued function.
This assumption is used only in section 4, formula (4.13).
If it is discarded, the part of section 4, which is based on this
assumption, should be changed. For example, if one assumes that
$h$ is sufficiently small, then the denominator
of (4.12) does not vanish for any $j\in \nsp$ and lemma 4.1 remains valid
without the assumption about real-valuedness of $h$.
\pn
Suppose now that for some
 $h\in C([0,T])$ problem (1.4)--(1.5) has a unique
solution
\beqn
%(1.8)
u\in C^2
([0,T];H)\cap C([0,T];{\cal D}(A)).
\eeqn
The assumption  $h\in C([0,T])$ is used in section 2 (see equation (2.3))
to identify $h$. Uniqueness of the solution to the
identification problem (IP), formulated below formula 
(1.9), is proved for
$h\in C([0,T])$ in section 3.
\pn
Assume that the kernel $h$ is  unknown and the data
\beqn
%(1.9)
g(t) = (u(t),\phi),\qq 0\le t\le T,
\eeqn
are measured.
\pn
The IP (identification problem) we study is: {\it given the data $(A,\phi,u_0,u_1,f,g)$,
with $g\in C^2([0,T])$, find a pair $(u,h)$ satisfying (1.4),
(1.5), (1.9).}
\pn
We note also that the exact data satisfy the additional
conditions
\beqn
%(1.10)
g(0)=(u_0,\phi),\qq g'(0)=(u_1,\phi),\qq g''(0)=(f(0),\phi).
\eeqn
If the the function $g(t)$ is not considered as an exact datum,
that is, a function of the form (1.9) where $u$ solves
(1.4)-(1.5) then
conditions (1.10)  are necessary for a solution to
problem (1.4), (1.5), (1.9) to exist.
\pn
We notice that the inverse problem we are going to study differs from the ones studied
in [9] and, e.g., [1], [2], and the very approach to our identification problem is new.
\pn
Consider now the more general equation depending on the {\it negative} parameter
$\mu$:
\beqn
%(1.11)
&& u''(t) = \mu Au(t) + \int_0^t h(t-s)Au(s)\,ds + f(t),\qq 0\le t\le T.
\eeqn
Observe that our identification problem can be viewed as the {\it limit} case
of the identification problem (1.11), (1.5), (1.9) as $\mu \to 0-$. We recall
that such problems have been widely studied under the assumption that $A$ is
selfadjoint and positive definite, cf., e.g, [1]--[8].

The main result, proved in section 2, is:

\theorem{1.1} If (1.1) and (1.3) hold, then the IP has at most
one solution.
\endproc
\med
\pn
An algorithm for finding $h$ and $u$ from the data is described 
in section 2.
\med

We can deal also with the first-order identification problem
consisting in determining a pair of functions $u:[0,T]\to H$ and
$l:[0,T]\to \rsp$ satisfying the Cauchy problem
\beqn
%(1.12)
u'(t) &=& l*Au(t)+f(t),\qq 0 \leq t \leq T,\\[2mm]
%(1.13)
u(0) &=& u_0.
\eeqn
as well as the additional conditions
\beqn
%(1.14)
g(t)=(u(t),\varphi),\qq 0\le t\le T,\qq l(0)\le 0,
\eeqn
where
\beqn
%(1.15)
f\in C^2([0,T];H),\q u_0\in {\cal D}(A),\q g(0)=(u_0,\phi)\neq
0,\q g\in C^3([0,T]).
\eeqn
Moreover, we assume that our data satisfy the additional conditions
\beqn
%(1.16)
g(0)=(u_0,\phi),\qq g'(0)=(f(0),\phi).
\eeqn
Such conditions are necessary for a solution to
problem (1.12)--(1.14) to exist. For exact data they are
satisfied automatically.
\pn
Differentiating both sides in (1.12) and taking (1.13) 
into account, it is immediate
to deduce that problem (1.12)--(1.14) 
is equivalent to problem (1.11), (1.5), (1.9)
with $\mu=l(0)$, $h=l'$, $f$ being replaced 
with $f'$ and $u_1=f(0)$. Since the case
corresponding to $l(0)<0$ has already been studied
in the literature, 
as was mentioned above, we can restrict ourselves 
to the study of the case $l(0)=0$. 
So conditions (1.14) can be replaced with the more
specific one
\beqn
%(1.17)
g(t)=(u(t),\varphi),\qq 0\le t\le T,\qq l(0)=0.
\eeqn
\pn
Finally we observe that a problem of 
the same type as (1.4), (1.5), (1.9) can
be treated similarly for the more general equation:
\beqn
%(1.18)
&& u''(t) = A_0u(t) + \int_0^t h(t-s)Au(s)\,ds + f(t),\qq 0\le t\le T.
\eeqn
Here $A_0:{\cal D}(A_0)\to H$ is a linear closed operator such that
\beqn
%(1.19)
{\cal D}(A_0)\supset {\cal D}(A),\qq \phi \in {\cal D}(A^*_0).
\eeqn
Further assume that the additional information
\beqn
%(1.20)
g_0(t) = (u(t),A_0^*\phi),\qq 0\le t\le T,
\eeqn
is available. We can consider the identification problem $IP_0$ related to
equations (1.18), (1.5), (1.9), (1.20) and 
to the data $(A,A_0,\phi,u_0,u_1,f,g,g_0)$.
Under the assumptions similar to those of theorem 1.1, one can uniquely
and algorithmically recover functions $h$ and $u$ from the data. However, the
existence of $u$ cannot be guaranteed, since the identification problem (1.18),
(1.5), (1.9), (1.20) is, in general, {\it overdetermined}. Yet, the existence
can be proved if, e.g., $\phi$ is a common eigenvector to $A$ and $A_0$.
\pn
We now describe the plan of the paper:
\pn
Section 2 is devoted to the existence and uniqueness of the unknown kernel $h$.
\pn
In section 3 the existence and uniqueness of the solution to the direct problem
(1.4), (1.5), with general closed selfadjoint operators satisfying (1.1)--(1.3),
are proved under suitable assumptions on the data.
\pn
In section 4 a mixed initial and boundary value problem is posed for the
operator equation (1.4) under the assumption that $A$ is a (closed) selfadjoint
and positive definite operator. Such a problem is solved under the assumptions
that the Fourier coefficients of the data decay sufficiently fast.
The results so found are then applied to the first-order equation $u'-l*Au=f$
with $l(0)=0$.
\pn
In section 5 some applications to linear integro-partial equations are considered.

\section{Uniqueness of the solution to IP}
\setcounter{equation}{0}
In this section we prove Theorem 1.1.

\proof  Multiply  both sides of (1.4) by $\phi$  and use properties (1.1)
and (1.3) to get:
\beqn
%(2.1)
 g''(t)- \l_0\int_0^t  h(t-s)g(s)\,ds = \psi(t),
 \q \psi(t) = (f(t),\phi),\q 0\le t\le T.
\eeqn
Equation (2.1) can be written as a linear
 Volterra integral equation for $h$:
\beqn
%(2.2)
&& \int_0^t h(t-s)g(s)\,ds
= \frac{g''(t)-\psi(t)}{\l_0}:=p(t),\qq 0\le t\le T.
\eeqn
Differentiate (2.2) and use (1.6) and (1.7) to get
\beqn
%(2.3)
&& g(0)h(t) + \int_0^t g'(t-s)h(s)\,ds = p'(t),\qq 0\le t\le T.
\eeqn
Note that the left-hand side of (2.2) is differentiable
even when $h\in C([0,T]),$ because  $g\in C^2([0,T]).$
Also, equation (2.2) shows that $p(t)\in C^1([0,T]),$
even  when $g\in C^2([0,T];H),$ because the left-hand side of
(2.2) and $\psi (t)$ are differentiable. Therefore, assuming
that $g\in C^2([0,T];H),$ and $f\in C^1([0,T];H),$ implies,
via equation (2.2), that $g\in C^3([0,T];H).$
 
Since $g(0)\neq 0$ by assumptions (1.7) and (1.10), equation (2.3)
is a Volterra second-kind equation with continuous kernel $g'$.
Therefore it is uniquely solvable 
in $C([0,T])$, since, according to our assumptions on the data,
$p'\in C([0,T])$.

It is well-known, that the solution to 
(2.3) can be obtained by iterations,
or, in analytic form, by the Laplace transform 
if one assumes $T=+\infty$.
It is also well-known that for any Volterra operator $V$ we have
$(cI+V)^{-1}=c^{-1}I+V_1$ if $c={\rm const}\neq 0$, $V_1$ 
being, in turn, a Volterra
operator.
\pn
We have assumed that the solution to (1.4), (1.5), with a known kernel
$h\in C([0,T])$, does exist and is unique. Therefore, if $h$ is found
(from (2.3)), then $u$ is uniquely found 
from (1.4), (1.5). Thus, theorem
1.1 is proved.\ \endproof
\med
\pn
An algorithm for the recovery of $h$ and $u$ 
from the data consists of solving (2.3)
for $h$ and then, once $h$ is known, solving (1.4), (1.5) for $u$.
\pn
Existence and uniqueness of the solution to 
(1.4), (1.5) is studied in section 3.
\pn
Finally, let us discuss problem 
$IP_0$ related to equations (1.9), (1.5), (1.18).
Multiply (1.18) by $\phi$ and get
\beqn
%(2.4)
&& g''(t) - g_0(t) - \l_0 \int_0^t h(t-s)g(s)\,ds = \psi (t),\qq 0\le t\le T,
\eeqn
where $\psi$ is defined in (2.1).
This equation can be reduced to equation (2.3) with function $q'$ replacing $p'$,
where
\beqn
%(2.5)
&& q(t) = \frac{g''(t)-g_0(t)-\psi(t)}{\l_0},\qq 0\le t\le T.
\eeqn
Thus $h$ is uniquely determined from the data for $IP_0$.
If one assumes that $\varphi$ is also an eigenvector of $A_0^*$ with an
eigenvalue $\lambda_{0,0}$, then the function $g_0$ in
formulas (2.4), (2.5) can be replaced by $\lambda_{0,0} g$.

A different approach to a study of (2.1) is given in section 4: rewrite
first (2.1) as
\beqn
%(2.6)
[g(t)-g(0) -tg'(0) - \psi^{(1)}(t)]/\lambda_0 =g*h^{(1)}(t),\qq 0\le t\le T,
\eeqn
where the superscript $(1)$ stands for convolution of a function $r$ with $t$, i.e.
\beqn
%(2.7)
r^{(1)}(t)=\int_0^t (t-s)r(s)\,ds,\qq 0\le t\le T.
\eeqn
Then solve the Volterra equation of the first kind (2.6) for $h^{(1)}$,
as it has been done above. If $h^{(1)}$ is found, then $h=(h^{(1)})''$.

To conclude this section we deal with problem (1.12)--(1.14). 
As before one
proves that $l\in C^1([0,T])$ is uniquely 
defined from the data. Indeed, from
(1.12) and (1.14) one derives  that
\beqn
%(2.8)
g'(t)=\lambda_0 l*g(t) +\psi(t),\qq \psi(t) = (f(t),\phi),
\qq 0\le t\le T,
\eeqn
where $\lambda_0$ is defined in (1.1)-(1.2). The case
of complex $\lambda_0 \in \csp \bs \{0\}$ can be treated similarly, as
explained in section 1.

From (2.8) one gets:
\beqn
%(2.9)
g*l(t) = [g'(t)-\psi(t)]/\lambda_0:=w(t), \qquad 0 \leq t \leq T.
\eeqn
Moreover, since (2.9) is a 
Volterra equation of the first kind for the unknown $l\in
C^1([0,T])$, and
it has at most one solution. Since $g(0)\neq 0$, and $w\in C^1([0,T]),$ then
differentiating (2.9) one gets the second kind Volterra equation for $l$.
\beqn
%(2.10)
g(0)l(t) + g'*l(t) = w'(t), \qquad 0 \leq t \leq T.
\eeqn
This yields the existence and uniqueness of $l$ and an algorithm for
the recovery of $l$, since the second kind Volterra equation can be
solved by iterations.

\remark{2.1} From formulas (2.10), (2.9) and (2.8) we easily compute the
initial value of $l$:
\beqn
%(2.11)
l(0)=\frac{g''(0)-(f'(0),\phi)}{g(0)\lambda_0}.
\eeqn
Hence, necessary conditions for equation (2.9) to admit a solution satisfying
$l(0)=0$ are:
\beqn
%(2.12)
g'(0)=(f(0),\varphi),\qq g''(0)-(f'(0),\phi)=0.
\eeqn

\section {Existence and uniqueness of the solution\\ to the direct problem (1.4), (1.5)}
\setcounter{equation}{0}
Let us assume that $h\in C
([0,T])$ is {\it known}, $A$ is selfadjoint and $E_\l$
is its resolution of the identity. The subspace $H_\L=E_\L$ is invariant with
respect to $A$ and $\|A\|_{H_\L}\le \L$.
\pn
Assume that the following hypothesis holds:
\beit
\item
[{\it H1}]\q $u_0,u_1\in H_\L$,\q $f\in C([0,T];H_\L)$.
\eeit
Applying $E_\L$ to (1.4), using {\it H1} and denoting $u_\L(t)=E_\L u(t)$,
$t\in [0,T]$, one gets equations (1.4) and (1.5) for $u_\L$. Problem (1.4),
(1.5) in $H_\L$ with a bounded operator $A$, $\|A\|_{H_\L}\le \L$ in $H_\L$,
is easily seen to be uniquely solvable, so existence and uniqueness of $u_\L$
follow. If the hypothesis {\it H1} holds with some $\L$, then it holds for
$H_\mu$ with any $\mu > \L,$ because $H_{\Lambda} \subset H_{\mu}$ for
$\mu> \Lambda$. Therefore existence and uniqueness of the solution $u$ to (1.4),
(1.5) is proved in any $H_\mu$, $\mu > \L$, provided that the hypothesis
{\it H1} holds.
\pn
This implies uniqueness of the solution to (1.4), (1.5) in $H$ if {\it H1} holds.
Indeed, assuming there are two solutions $u_1$ and $u_2$ to problem (1.4), (1.5),
one concludes from the above argument that $\|u_1-u_2\|_{H_\mu}=0$ for all
$\mu \ge \L$. Since $0=\lim_{\mu\to +\infty}\,\|u_1-u_2\|_{H_\mu}=\|u_1-u_2\|_H$,
it follows that $u_1=u_2$.

Our argument proves that {\it the homogeneous direct problem (1.4), (1.5) has
only the trivial solution}.

Since  $f=0$, $u_0=0$ and $u_1=0$ satisfy {\it H1}, problem (1.4), (1.5) with
a selfadjoint $A$ and $h\in C^2([0,T])$ has at most one solution in $H$.
Indeed, if it has two solutions, their difference, $u$, solves the homogeneous
problem (1.4), (1.5). Consequently $u(t)=\int_0^t h^{(1)}(t-s)Au(s)\,ds$,
where $h^{(1)}(t)$ is defined by formula (2.7) with $r$ replaced by $h$.

Fix now an arbitrary $\Lambda < \infty$ and apply $E_{\Lambda}$ to get
$u_{\Lambda}(t)=\int_0^t h^{(1)} (t-s)A_{\Lambda}u_{\Lambda}(s)\,ds$,
where $A_{\Lambda}:=E_{\Lambda}A=E_{\Lambda}AE_{\Lambda},$
and we have used the formula $E_{\Lambda}A=AE_{\Lambda}$,
$E_{\Lambda}^2=E_{\Lambda}.$
Since $A_{\Lambda}$ is a bounded linear operator, it follows that
$u_{\Lambda}=0.$ Since $\Lambda$ is arbitrary, this implies $u=0$.

Existence of the solution requires special assumption on $f$, $u_0$ and $u_1$.
Since usually $f$, $u_0$, $u_1$ are at our disposal when we study the inverse problem,
assumption {\it H1} is not restrictive and is quite natural: if $A$ is known, then
$E_\L$ and $H_\L$ are known, and one can choose the data $f$, $u_0$, $u_1$ in
$H_\L$. Moreover if the data are noisy, that is $f$, $u_0$, $u_1$ are known up to a
(known) error $\dl$, i.e
\beqn
%(3.1)
\| u_{0,\dl}-u_0\|\le \dl,\qq \|u_{1,\dl}-u_1\|\le \dl,
\qq \|f_{\dl}-f\|_{C([0,T];H)}\le \dl,
\eeqn
then one can use the data $E_\L f_{\dl}$, $E_\L  
u_{0,\dl}$, $E_\L u_{0,\dl}$
which satisfy {\it H1}. Since $E_\L$ is known, 
computation of $E_\L f_{\dl}$,
$E_\L u_{0,\dl}$, $E_\L u_{0,\dl}$ presents no difficulties. 
In these arguments
we assume that $A$ is given exactly.
\pn
If one wants to weaken assumption {\it H1}, one can allow the data to have a non-zero
component in $H\ominus H_\L$, but this component must have coefficients exponentially
decaying as $\l \to +\infty$.

Let us summarize the results of this section:

\theorem{3.1} Let $A=A^*$ be a possibly unbounded operator and let
$h\in C([0,T])$. Then problem (1.4), (1.5) has at most one solution in
$C^2([0,T];H)\cap C([0,T];{\cal D}(A))$. If in addition the hypothesis H1
holds,  then problem (1.4), (1.5) has a solution in $C^2([0,T];H_{\Lambda})$
and this solution is unique.
\endproc
\med
\pn
Recall now that, if $f\in C^1([0,T];H)$, and
$l\in C^1([0,T])$ with $l(0)=0$, then 
the direct problem (1.12), (1.13)
with $l(0)=0$, is equivalent to the second-order Cauchy problem
\beqn
%(3.2)
&& u''(t) = \int_0^t h(t-s)Au(s)\,ds + f'(t),\qq 0\le t\le T,\\[2mm]
%(3.3)
&&  u(0)=u_0,\q u'(0)=f(0),
\eeqn
where $h(t)=l'(t)$. Further assume
\beit
\item
[{\it H2}]\q $u_0\in H_\L$,\q $f\in C^1([0,T];H_\L)$.
\eeit
Then from (3.2), (3.3) and Theorem 3.1 we get the following
theorem:

\theorem{3.2} Let $A=A^*$ be a possibly unbounded operator and let
$l\in C^1([0,T]),\q l(0)=0.$ Then problem (1.12), (1.13) has at
most one solution in
$C^2([0,T];H)\cap C([0,T];{\cal D}(A))$. If in addition the hypothesis H2 holds,
then problem (1.12), (1.13) has a solution in $C^2([0,T];H_{\Lambda})$ and this
solution is unique.
\endproc
\med
\pn

\section{A mixed problem for equation (1.4)}
\setcounter{equation}{0}
In this section the solution to (1.4) which satisfies the boundary conditions
\beqn
%(4.1)
u(0)=u_0,\qq u(T)=u_2,
\eeqn
is studied. Note that $u''(0)=f(0),$ as follows from (1.4)
if $u\in C^2([0,T];H)$.

We will show that, under suitable assumptions, the data
\beqn
%(4.2)
(A,\phi,u_0,u_2,f,g)
\eeqn
determine the pair $(u,h)$ uniquely.
\pn
Let us assume that the operator $A$ does not depend on time, $A=A^*$
and $\{\phi_j\}_{j=1}^{+\infty}$ is an orthonormal basis of $H$ such that
$A\phi_j=\l_j\phi_j$, $j\in \nsp$, $\{\lambda_j\}_{j=1}^{+\infty}$
being a positive nondecreasing sequence diverging to $+\infty$

 If (1.4) is solvable in $C^2([0,T];H)$, then
\beqn
%(4.3)
u(t) = \sum_{j=1}^{+\infty}\,{\what u}_j(t)\phi_j,\qq {\what u}_j(t)=(u(t),\phi_j),
\qq 0\le t\le T.
\eeqn
Hence the Fourier coefficients ${\what u}_j$ solve the {\it scalar} boundary
value problems
\beqn
%(4.4)
&& {\what u}_j''(t)-\l_jh*{\what u}_j(t)={\what f}_j(t),\\[2mm]
%(4.5)
&& {\what u}_j(0)={\what u}_{0,j}:=(u_0,\phi_j),\qq {\what
u}_j(T)={\what u}_{2,j}:=(u_2,\phi_j),
\eeqn
where $f*g(t)=\int_0^t f(t-s)g(s)ds$.
\pn
Let
\beqn
%(4.6)
{\what u}_j'(0)=c_j,\qq {\wtil f}_j(t) = \int_0^t (t-s){\what f}_j(s)\,ds.
\eeqn
Then (4.4) implies
\beqn
%(4.7)
{\what u}_j(t) = {\what u}_{0,j} + tc_j + {\wtil f}_j(t) + \l_jh^{(1)}*{\what u}_j(t),
\qq 0\le t\le T,
\eeqn
where
\beqn
%(4.8)
h^{(1)}(t) = \int_0^t (t-s)h(s)\,ds,\qq 0\le t\le T.
\eeqn
Define the Volterra operators $K_j$, $j\in \nsp$, by the formulas
\beqn
%(4.9)
(I-\l_jH)^{-1}=I+\l_jK_j,\qq Hf:=h^{(1)}*f,\qq K_jf:=k_j*f,\qq \forall
j\in \nsp.
\eeqn
Note that $h\in C([0,T])$ and $f\in C([0,T];H)$ imply $k_j\in C([0,T])$
for all $j\in \rsp$.
\pn
Then (4.7) and (4.9) imply, for any $t\in [0,T]$,
\beqn
%(4.10)
{\what u}_j(t) = {\what u}_{0,j} + tc_j + {\wtil f}_j(t) + \l_j(k_j*1)(t){\what u}_{0,j} +
\l_j(k_j*t)(t)c_j + \l_jk_j*{\wtil f}_j(t).
\eeqn
To satisfy condition ${\what u}_j(T)={\what u}_{2,j}$ for any $j\in \nsp$ it
is sufficient to assume
\beqn
%(4.11)
\Big|T + \l_j\int_0^T k_j(T-s)s\,ds\Big| > 0,\qq \forall j\in \nsp.
\eeqn
Under such an assumption from (4.10) it easily follows
\beqn
%(4.12)
c_j = {{\what u}_{2,j} - {\what u}_{0,j} - {\wtil f}_j(T) - \l_j\int_0^T k_j(T-s)[{\what u}_{0,j}
+ {\wtil f}_j(s)]\,ds \over T + \l_j\int_0^T k_j(T-s)s\,ds}.
\eeqn
If
\beqn
%(4.13)
h(t)\ge 0, \qq 0\le t\le T, \qq \l_j>0,\qq \forall j\in \nsp,
\eeqn
then (4.9) shows that
\beqn
%(4.14)
K_j=\sum_{m=1}^{+\infty}\, \l_j^{m-1}H^m \q \Longleftrightarrow \q
k_j=\sum_{m=1}^{+\infty}\, \l_j^{m-1}(h^{(1)}*)^{m-1}h^{(1)},
\eeqn
so $k_j(t)\ge 0$ for all $t\in  [0,T]$ and all $j\in \nsp$. Therefore condition (4.11)
is satisfied.
\pn
Therefore we have proved

\lemma{4.1} If (4.13) holds, then, for any $j\in \nsp$, problem (4.4), (4.5) is
solvable for any triplet $({\what u}_{0,j},{\what u}_{2,j})\in \rsp^2$,
${\what f}_j\in C([0,T])$, and its solution is unique.
\endproc

\remark{4.1} According to equation (4.3) condition (4.11) is satisfied if we
assume, e.g., that the data fulfill the following inequalities
\beqn
%(4.15)
g(0)g'(t)< 0,\qq \l_0g(0)[g'''(t)- \l_0 g'(t)-(f'(t),\phi)]
> 0,\qq 0\le t\le T.
\eeqn
Indeed, it suffices to show that (4.15) implies $k_j(t)\ge 0$ for any $t\in
[0,T]$. The  first of the conditions in (4.15) and the
definition
of function $p$ in (2.2) imply that the solution $h$ to equation (2.3), rewritten
as a fixed-point equation, is non-negative in $[0,T]$. Hence, the same holds for
$h^{(1)}$ according to (4.8). As a consequence, since $\l_j > 0$ for all
$j\in \nsp$, from (4.14) we immediately deduce our assertion.
\med

In general, the solution to (4.4) may not exist if the denominator of (4.12)
vanishes for some integer $j_0$. For the solution to (4.4) to exist in this
case it is necessary and sufficient that the numerator of (4.12) vanishes also.
If this is the case, the solution exists but is not unique: it is of the form
(4.10) with $c_{j_0}$ arbitrary.
\pn
For the series (4.3), representing the solution to (1.4) and (4.1),  to converge
in $C^2([0,T];H)$ \newline  $\cap C([0,T];{\cal D}(A))$ it is necessary and sufficient that
\beqn
%(4.16)
\sup_{0\le t\le T}\,\sum_{j=1}^{+\infty}\,
\big\{|{\what u}_j''(t)|^2+\l_j^2|{\what u}_j(t)|^2\big\} < +\infty,
\eeqn
where ${\what u}_j$ is defined by (4.10) and (4.12).
\pn
Condition (4.16) is equivalent to requiring
\beqn
%(4.17)
\sup_{0\le t\le T}\,\sum_{j=1}^{+\infty}\,\l_j^2|{\what u}_j(t)|^2 < +\infty,
\eeqn
as follows from (4.4).
\pn
Condition (4.17) is trivially satisfied if, for example,
\beqn
%(4.18)
{\what u}_j(T)={\what u}_{0,j}={\wtil f}_j(s)=0,\q \mbox{for } j\ge J,\ 0\le s \le T,
\eeqn
where $J$ is an arbitrarily large, fixed integer.
\pn
Condition (4.18) is sufficient but not necessary for (4.17) to hold.
\pn
Let us derive a less restrictive sufficient condition for (4.17) to hold,
which is close to a necessary one. Denote
\beqn
%(4.19)
|{\what u}_{0,j}|^2 +T^2 |c_j|^2 + \sup_{0\le t\le T}\, |{\wtil f}_j(t)|^2
:= L_j,\qq \forall j\in \nsp.
\eeqn
Formula (4.10) implies
\beqn
%(4.20)
|{\what u}_j(t)|^2 \le 6\Big(1 + T\l_j^2 \int_0^T |k_j(t)|^2\,dt\Big)
L_j, \qq \forall j\in \nsp,
\eeqn
where the Cauchy inequality and the following elementary inequalities
were used:
$$
\Big(\sum_{j=1}^n a_j\Big)^2\leq n \Big(\sum_{j=1}^n a_j^2\Big),\qq a_j\geq 0,
$$
with $n=2$ and $n=3$.
\pn
Let us now estimate $k_j$. If we denote
\beqn
%(4.21)
T\int_0^T |h(t)|\,dt := M,
\eeqn
then
\beqn
%(4.22)
\sup_{0\le t\le T}\, h^{(1)}(t)\le M.
\eeqn
If we set
\beqn
%(4.23)
h_m(t) = \int_0^t h^{(1)}(t-s)h_{m-1}(s)\,ds,\qq h_1(t)=
h^{(1)}(t),\qq 0\le t\le T,
\eeqn
where $h^{(1)}(t)$ is defined in (4.8), then, by induction, one gets
\beqn
%(4.24)
h_m(t) \le  M\frac{(Mt)^{m-1}}{(m-1)!},\qq 0\le t\le T,\ \forall m\in \nsp.
\eeqn
Therefore (4.14) implies
\beqn
%(4.25)
0\le k_j(t)\le M\exp(\l_jMt),\qq 0\le t\le T,\ \forall j\in \nsp,
\eeqn
and
\beqn
%(4.26)
\int_0^T |k_j(t)|^2\,dt \le M^2\frac{\exp(2\l_jMT)-1}{2\l_jM} \le
M\frac{\exp(2\l_jMT)}{2\l_j},\ \forall j\in \nsp.
\eeqn
Since the positive nondecreasing 
sequence $\{\lambda_j\}_{j=1}^{+\infty}$
diverges to $+\infty$, formulas (4.20) and (4.26) imply
\beqn
%(4.27)
\sup_{0\le t\le T}\,|{\what u}_j(t)|^2\le
6 L_j \big[1+0.5 MT\l_j\exp(2\l_jMT)\big].
\eeqn
From (4.27) it follows that (4.17) holds if
\beqn
%(4.28)
\sum_{j=1}^{+\infty}\,L_j\l_j^3\exp(2\l_jMT)<+\infty.
\eeqn
Condition (4.28) means, 
roughly speaking, that for the series (4.28) to
converge the sequences $\{{\what u}_{0,j}\}_{j=1}^{+\infty}$,
$\{{\what u}_{2,j}\}_{j=1}^{+\infty}$ and 
$\{{\wtil f}_j\}_{j=1}^{+\infty}$
related to the Fourier coefficients of the data must decay 
faster than
$\{\l_j^{-3/2}\exp(-\l_jMT)\}$ as $j\to \infty.$.
\pn
Condition (4.28) is not far from a necessary condition since the above estimates
were not too crude.

Let us summarize the result:
\pn
\theorem{4.1}  If (4.13) and (4.28) hold, then, for any non-negative
$h\in C([0,T])$ the solution of the problem (1.4) and (4.1) in
$C^2([0,T];H)\cap C([0,T];{\cal D}(A))$ does exist and is unique.
\endproc
\med

To conclude this section we consider the identification problem consisting of
recovering the pair $u:[0,T]\to H$ and $l:[0,T]\to \rsp$ satisfying the first-order
Cauchy problem
\beqn
%(4.29)
u'(t) &=& l*Au(t)+f(t),\qq 0 \leq t \leq T,\\[2mm]
%(4.30)
u(T) &=& u_2.
\eeqn
as well as the extra data
\beqn
%(4.31)
g(t)=(u(t),\varphi),\qq l(0)=0.
\eeqn
The assumptions about $A$ are the same as at the 
beginning of this section, and $f\in C^1([0,T];H)\q
g\in C^3([0,T]).$
\pn
Moreover, recalling (2.12), we 
assume that the data satisfy the conditions
\beqn
%(4.32)
g'(0)=(f(0),\varphi),\q g''(0)-(f'(0),\phi)=0,
\q g(T)=(u_2,\varphi),\q
g(0)\neq 0.
\eeqn

Recall now that, since $l(0)=0$, problem (4.29), (4.30) is equivalent to
the second-order Cauchy problem
\beqn
%(4.33)
&& u''(t) = \int_0^t h(t-s)Au(s)\,ds + f'(t),\qq 0\le t\le T,\\[2mm]
%(4.34)
&& u'(0)=f(0),\q u(T)=u_2,
\eeqn
where $h(t)=l'(t)$.
\pn
Observe that the present problem differs from the one just studied
only by the initial condition: $u'(0)=f(0)$ replaces
$u(0)=u_0$.
Moreover, since we have no initial condition $u(0)=u_0$, we need the explicit
requirement $g(0)\neq 0$ (cf. (4.32)) to ensure that equation (2.10) is
actually of the second kind.
\pn
Reasoning as at the beginning of this section, we easily get that formula
(4.7) is now replaced by
\beqn
%(4.35)
{\what u}_j(t) &=& c_j + t{\what f}_j(0) + \int_0^t (t-s){\what f}'_j(s)\,ds
+\l_jh^{(1)}*{\what u}_j(t)\no \\[2mm]
&=& c_j + \int_0^t {\what f}_j(s)\,ds + \l_jh^{(1)}*{\what u}_j(t),
\eeqn
where $c_j$ stands for the unknown value ${\what u}_j(0)$ and, in our case,
\beqn
%(4.36)
h^{(1)}(t)= \int_0^t (t-s)l'(s)\,ds = \int_0^t l(s)\,ds,\qq 0\le t\le T.
\eeqn
Introducing the kernels $k_j$ defined in (4.9) finally we get
the representation formulas
\beqn
%(4.37)
{\what u}_j(t) = c_j + 1*{\what f}_j(t) + c_j\l_j(k_j*1)(t)
+ \l_jk_j*(1*{\what f}_j(t)).
\eeqn
Moreover, the solvability condition (4.11) changes to
\beqn
%(4.38)
\Big|1 + \l_j\int_0^T k_j(s)\,ds\Big| > 0, \qq \forall j\in \nsp.
\eeqn
Using (4.38), we obtain
\beqn
%(4.39)
c_j = {{\what u}_{2,j} - \int_0^T {\wtil f}_j(s) 
\Big[ 1 + \l_j\int_0^{T-s} k_j(\s)\,d\s \Big]ds
\over 1 + \l_j\int_0^T k_j(s)\,ds},\qq \forall j\in \nsp.
\eeqn

\remark{4.2} According to equation (2.7) condition (4.38) is satisfied if we
assume, e.g., that the data fulfill the following inequalities
\beqn
%(4.40)
g(0)g'(t)< 0,\qq \l_0g(0)[g''(t)-(f'(t),\phi)]> 0,\qq 0\le
t\le T.
\eeqn
Indeed, (4.40) implies that the solution $l$ to equation (2.10) is
non-negative
in $[0,T]$. Then (4.36) implies $h^{(1)}(t)\ge 0$ for all $t\in [0,T]$.
Finally, from (4.14) we deduce that $k_j$ is non-negative in $[0,T]$ for any
$j\in \nsp$.
\med

To summarize our basic result for first-order integro-differential
equations we need the notation
\beqn
%(4.41)
|{\what u}_{2,j}|^2 + |c_j|^2 + \sup_{0\le t\le T}\, |{\wtil f}_j(t)|^2 := L_j,
\qq \forall j\in \nsp.
\eeqn

\theorem{4.2} If $f\in C^1([0,T])$ and (4.28), (4.32) hold, then, for any
non-negative $l\in C^1([0,T])$ the solution $u$ of the problem (4.29),
(4.30) in $C^2([0,T];H)\cap C([0,T];{\cal D}(A))$ does exist and is unique.
\endproc
\med
\pn
\section{Applications}
\setcounter{equation}{0}

{\bf Example 1.} We apply the previous abstract result for second-order
integro-differential equations to the identification problem:
{\it determine two functions $u:[0,T]\times \Om\to \rsp$ and $h:[0,T]\to \rsp$
satisfying the following equations for some $m\in \{0,1\}$:}
\beqn
%(5.1)
\hskip -2truecm
&& D_t^2u(t,x) = \int_0^t h(t-s){\cal A}(x,D_x)u(s,x)\,ds
\nonumber \\[2mm]
&& \hskip 2.2truecm + f(t,x),\qq (t,x)\in [0,T]\times \Om,\\[2mm]
%(5.2)
\hskip -2truecm
&& u(mT,x)=u_{2m}(x),\q D_tu(0,x)=u_1(x),\qq x\in \Om,\\[2mm]
%(5.3)
\hskip -2truecm
&& {\cal B}u(t,x)=0,\qq (t,x)\in (0,T)\times \partial \Om,\\[2mm]
%(5.4)
\hskip -2truecm
&& \int_\Om \phi(x)u(t,x)\,dx=g(t),\qq 0\le t\le T.
\eeqn
Here $\Om$ denotes a bounded domain in $\rsp^n$ with boundary
$\partial \Om$ of class $C^{1,1}$, while
\beqn
%(5.5)
&& {\cal A}(x,D_x)u= - \sum_{i,j=1}^n\,
D_{x_i}(a_{i,j}(x)D_{x_j}u)+a_{0,0}(x)u,\\[2mm]
%(5.6)
&& {\cal A}(x,D_x)\phi(x)=\l_0\phi(x),\q \l_0\neq 0.
\eeqn
Assume that:
\beqn
%(5.7)
\hskip -2truecm
&& \overline {a_{i,j}(x)}=a_{j,i}(x)=a_{i,j}(x)\in C^{0,1}
({\ov \Om}),\q
i,j=1,\ldots,n,\q a_{0,0}\in L^\infty(\Om),\\[2mm]
%(5.8)
\hskip -2truecm
&&\sum_{i,j=1}^n\, a_{i,j}(x)\xi_i\xi_j\ge \mu \sum_{j=1}^n |\xi|^2 >0,
\qq \forall x\in {\ov \Om},\ \forall \xi \in \rsp^n.
\eeqn
Choose, for example,
\beqn
%(5.9)
{\cal B}u(x)=u(x),\qq x\in \partial \Om,
\eeqn
\pn
and
\beqn
%(5.10)
{\cal D}(A) = \{u\in H^2(\Om): u=0\ \mbox{on } \partial \Om\},\q
Au(x)={\cal A}(x,D_x)u(x),\q u\in {\cal D}(A).
\eeqn

\remark{5.1} Equation (5.1) is not of standard type, elliptic or hyperbolic.
Also when the kernel $h$ has a fixed sign, the behavior
 of the solution
can be very wild (cf. section 4). 
At any rate, the behavior strongly depends
on the type of the prescribed conditions: initial, 
boundary or mixed ones.
\med
\pn
The results of sections 1--4 are applicable to problems (5.1)--(5.4) with
$m=0$ or $m=1$ and $H=L^2(\Om)$. Such results ensure, under explicit
conditions on the data, the existence and the uniqueness of a solution
$(u,h)$ to (5.1)--(5.4) and give an algorithm for recovery of $h$ from the
data.
\pn
Note that the existence of the non-zero eigenvalue of $A^*$
follows from the selfadjointness of $A$ and the known results
about the eigenvalues of elliptic operators.
\med
\pn
{\bf Example 2.} We apply the previous abstract result for first-order
integro-differential equations to the following identification problem:
{\it determine two functions $u:[0,T]\times \Om\to \rsp$ and
$l:[0,T]\to \rsp$ satisfying the following equations for some (fixed)
$m\in \{0,1\})$:}
\beqn
%(5.11)
\hskip -2truecm
&& D_tu(t,x) = \int_0^t l(t-s){\cal A}(x,D_x)u(s,x)\,ds
+ f(t,x),\qq (t,x)\in [0,T]\times \Om,\\[2mm]
%(5.12)
\hskip -2truecm
&& u(mT,x)={\what u}_{2m}(x),\qq x\in \Om,\\[2mm]
%(5.13)
\hskip -2truecm
&& {\cal B}u(t,x)=0,\qq (t,x)\in (0,T)\times \partial \Om,\\[2mm]
%(5.14)
\hskip -2truecm
&& \int_\Om \phi(x)u(t,x)\,dx=g(t),\qq 0\le t\le T,\qq l(0)=0.
\eeqn
Here $\Om$ and ${\cal A}(x,D_x)$ enjoy the same properties as in Example 1.
\pn
The results of sections 1--4 are applicable to problems (5.11)--(5.13) with
$m=0$ or $m=1$ and $H=L^2(\Om)$. Such results ensure, under explicit conditions
on the data, the existence and the uniqueness of a solution $(u,h)$ to (5.11)--(5.14)
and give an algorithm for recovery of $l$ from the data.

\end{document}